\documentclass[12pt]{article}

\usepackage[T2A]{fontenc}
\usepackage[cp1251]{inputenc}
\usepackage[english]{babel}

\usepackage{amscd}
\usepackage{amssymb}
\usepackage{amsmath}
\usepackage{graphicx}
\usepackage{graphics}
\usepackage{ifthen}
\usepackage{truncate}
\usepackage{amsopn}           
\usepackage{longtable}

\begin{document}
\title{On bifurcations in degenerate resonance zones}%
\author{A.D. Morozov \\
\small  Lobachevsky  State University of Nizhny Novgorod \\
\small Department of Mathematics and Mechanics\\
 \small 603950, Russia, Nizhny Novgorod, Gagarin Ave., 23\\
\small  E-mail: morozov@mm.unn.ru}%
\date{}%

\maketitle \noindent MSC2000 numbers: 34C15, 34C28, 37G35

\medskip


\begin{abstract}
For Hamitonian systems with 3/2 degrees of freedom close to
nonlinear integrable and for symplectic maps of the cylinder,
bifurcations in degenerate resonance zones are discussed\footnote
{The paper is based on the talk the author presented at the
International Conference ``Dynamics, Bifurcation, and Strange
Attractors'' dedicated to the memory of L.P.Shil'nikov (Nizhny
Novgorod, Russia, July 1-5, 2013). Only the stuff devoted to
degenerate resonances was included.}
\end{abstract}

\section{Introduction}
The study of orbit behavior in nonlinear nearly integrable
Hamiltonian systems with 3/2 degrees of freedom near the levels
close to a resonance one was heavily advanced by L.P. Shil'nikov
in the paper \cite{MorSh1983} (jointly with the author of this
paper). These investigations were continued in numerous papers
(references can be found in books \cite{GH}-\cite{Mor2005}). In
the papers mentioned resonances were classified as being passable,
partially passable and non-passable. Also they can be
non-degenerate or degenerate ones. Perturbations we consider can
be Hamiltonian or non-Hamiltonian. Degenerate resonances occur in
systems where the frequency of periodic orbits for an initial
one-degree-of-freedom unperturbed Hamiltonian system is a
non-monotone function of the value of Hamiltonian. The analysis of
degenerate resonances was initiated in \cite{MorSh1983}. After
that much work was done to analyze the degenerate resonances (see,
for instance, \cite{MorB1999}-\cite{Mor2008}, \cite{Mor2005}).

For symplectic maps on the cylinder with a non-monotone rotation number, degenerate
resonances were studied first in \cite{HowH1984}. Here, we address mostly
to bifurcations in degenerate resonance zones.

For systems under consideration, the motions in resonance zones
(we call by this term some region near the resonance
level of the unperturbed Hamiltonian) can be observed for
 a planar Poincar\'e map constructed by means of a
computer (using, for instance, the package WinSet \cite{MD}).
For sufficiently small perturbations, the
pictures obtained in this way are in a good agreement with the analysis of averaged
systems. In the case of Hamiltonian perturbations, the related Poincar\'e
map preserves of course the area.

As for area preserving maps of the cylinder is concerned,  the
Chirikov map (or the standard map, in the other terminology) is
the most popular. This is a map with a monotone rotation number.
Also, maps were studied where the rotation function is nonlinear
and non-monotone.

As the first example, maps with the quadratic rotation function
$f(u)=p_1u+p_2xu^2$ (or $f(u)=p_1+p_2u^2$) with parameters
$p_1,\,p_2$ were considered. Such the maps were called to be
``standard maps with a non-monotone rotation'' \cite{HowH1984}.
When studying, these maps were approximated by Hamiltonian flows
and bifurcations of reconnecting separatrices were examined
related with two scenarios: 1) the formation of ``loops'' and 2)
the formation of ``vortex pairs'' \cite{HowH1984}-\cite{Petrisor}.
Some related features were observed for the map itself as well.
However, no answer was found so far to the following question:
which resonances give rise to the formation of vortex pairs? For
instance, in \cite{HowHum1995} vortex pairs arose for resonances
with even number $2$  but not for resonances with other integers.

When studying Hamiltonian systems with 3/2 degrees of freedom
close  to nonlinear integrable, the second scenario has been
numerically realized \cite{Mor2008}. However, it was proved for
sufficiently small perturbations that ``vortex pairs'' are absent
there and only the first scenario for the reconnection of
separatrices is possible. The existence of vortex pairs in
\cite{Mor2008} is explained by the fact that the perturbation was
not ``sufficiently small''. In this paper, we study in details the
bifurcations in degenerate resonance zones, using the averaged
system of the second approximation. This study answers the
question on the formation of vortex pairs.

\section{Averaging}
Let us consider a system
\begin{equation}\label{8.1}
\begin{split}
\dot{x} &= {\partial H(x,y)\over \partial y} + \varepsilon g(x,y,\nu t)  \\
\dot{y} &= - {\partial H(x,y)\over \partial x} + \varepsilon
f(x,y,\nu t),
 \end{split}
\end{equation}
where $\varepsilon>0$ is a small parameter, functions
$H(x,y)$, $g(x,y,\nu t)$, $f(x,y,\nu t)$ are sufficiently smooth
(or analytic) in variables $x,\,y$ varying in some domain
$D\subset R^2$ (or $D\subset S^1\times R^1$) and continuous and
$2\pi/\nu$-periodic in $t$, $\nu$ is a
parameter (the perturbation frequency). We assume the unperturbed system to be nonlinear and to possess a cell $D_0$ filled with periodic  orbits. In $D_0$ the action-angle variables $I$,
$\theta$ are introduced. As a result one obtains the system
\begin{equation}\label{8.2}
\begin{split}
 \dot{I}& = \varepsilon [f(x,y,\varphi)x^{^\prime
}_{\theta }- g(x,y,\varphi )y^{^\prime }_{\theta
}]\equiv \varepsilon F(I,\theta ,\varphi )\qquad  \\
\dot{\theta } &= \omega (I)+\varepsilon [-f(x,y,\varphi
)x^{^\prime }_{I}+
 g(x,y,\varphi )y^{^\prime }_{I}]\equiv \omega (I)+\varepsilon G(I,\theta ,\varphi
 )\\
 \dot{\varphi }&=\nu.
 \end{split}
\end{equation}

Studying resonances play an essential role when analysing system
(\ref{8.2}). As is known, by a resonance in the system it  is
understood the presence of the following relation
\begin{equation}\label{res}
\omega (I)=(q/p)\nu ,
\end{equation}
where $p,\,q$ are co-prime integers. We denote the related value
$I$  in (\ref{res}) as $I_{pq}$.

If the following conditions hold
\begin{equation}\label{degres}
\omega^{(k)}(I_{0})=0,\, k=1,2,\dots,j-1; \, \omega^{(j)}(I_{0})
\neq 0, \quad j>1,
\end{equation}
we call the level $I=I_0$ to be {\it degenerate} and the number
$j$  is its degeneracy order. If $\omega^{\prime}(I_0)\neq 0$
(that is, $j=1$), then the level $I=I_0$ is called to be {\it
non-degenerate}.

By the resonance zone we mean henceforth a neighborhood
$$
U_{\varepsilon ^s}=\{(I,\theta ): I_{pq}-c\varepsilon ^s
<I<I_{pq}+c\varepsilon ^s , 0\le \theta <2\pi , c=\hbox{const}>0,
 s=1/(1+j)\}
$$
of an individual resonance level $I=I_{pq}$.  In this
neighborhood, according to \cite{Mor2005}, the behavior of
solutions to the initial system (\ref{8.2}) is described (up to
terms $O(\varepsilon^{1+s})$) by the autonomous system
\begin{equation}\label{8.27}
\begin{split}
\dot u &=\varepsilon ^{1-s}A_{0}(v,I_{pq})+
\varepsilon P_{0}(v,I_{pq})u  \\
 \dot v &=\varepsilon ^{1-s}b_{j}u^{j}+
 \varepsilon (b_{j+1}u^{j+1}+Q_{0}(v,I_{pq})),\, s=1/(1+j),
  \end{split}
\end{equation}
where $b_{j}=\omega ^{(j)}(I_{pq})/j!$,
\begin {equation} \label{eq9}
  A_0(v,I_{pq}) = {1\over 2\pi p} \int_0 ^ {2\pi p} F (I _ {pq},
  v +\frac{q}{p} \varphi, \varphi) d\varphi,
\end {equation}
\begin{equation}\label{eq91}
P_{0}(v,I_{pq})={1\over 2\pi p}\int^{2\pi p}_{0}[\partial
F(I_{pq},v+q\varphi /p, \varphi)/\partial I]d\varphi ,
\end{equation}
\begin{equation}\label{eq92}
Q_{0}(v,I_{pq})={1\over 2\pi p}\int^{2\pi
p}_{0}G(I_{pq},v+q\varphi /p,\varphi )d\varphi .
\end{equation}
Functions $A_0(v,I_{pq})$, $P_0(v,I_{pq})$, $Q_0(v,I_{pq})$
have a common period $2\pi/p$ in $v$. Let us write down
functions $A_0(v,I_{pq})$, $P_0(v,I_{pq})$ in the form
\begin{equation}\label{eq93}
\begin{split}
A_0(v,I_{pq}) &=\widetilde{A_0}(v,I_{pq}) +B_0(I_{pq}),\\
P_{0}(v,I_{pq})&=\widetilde{P_0}(v,I_{pq}) +B_1(I_{pq}),
\end{split}
\end{equation}
where $B_0(I_{pq})$ is the average value of $A_0(v,I_{pq})$ and
$B_1(I_{pq})$ is the same for $P_0(v,I_{pq})$.

Now consider the equation
\begin{equation}\label{eq95}
\widetilde{A_0}(v,I_{pq}) +B_0(I_{pq})=0,
\end{equation}
whose roots define coordinates $v$ of equilibria for the truncated
system. Following \cite{MorSh1983} a resonance level $I=I_{pq}$ is
called to be: i) passable, if the equation (\ref{eq95}) has not
real roots; ii) partially  passable if (\ref{eq95}) has only
simple real roots and $B_0(I_{pq})\neq 0$; iii) non-passable if
$B_0(I_{pq})= 0$.

System (\ref{8.27}) is usually called to be the averaged system of the
second approximation of the averaging method.

\section{Bifurcations in degenerate resonance zones}

Consider  first the case of Hamiltonian perturbations. Then the
initial system (\ref{8.1}) is Hamiltonian and the following
identities hold
   \begin{equation}\label{8.28}
P_{0}+Q^{^\prime }_{0v}\equiv  0, \quad B_0=B_1=0.
\end{equation}
Hence, system (\ref {8.27}) is also Hamiltonian. In this case,
the first approximation system plays the main role in the analysis
of phase curves of system (\ref {8.27})
\begin{equation}\label{eq94}
\begin{split}
&\dot{u} = \varepsilon ^{1-s}A_0(v,I_{pq})  \\
& \dot{v} =  \varepsilon ^{1-s}b_{j} u^{j}б \qquad j\geq 2.
\end{split}
\end{equation}

The explicit calculation of $A(v,I_{pq})$ for main examples of nonlinear
systems is relied on the
problem of finding unperturbed solution $x(\theta,I_{pq}),
y(\theta,I_{pq}) $ which, in turn, is intimately connected with the
problem of inversion for hyper-elliptic integrals.

System (\ref{eq94}) has only degenerate equilibria.  Therefore, along
with system (\ref{eq94}) its deformation is investigated
 \begin{equation}\label{eq12}
\begin {split}
  \dot u &= \varepsilon ^{1-s} A (v, I _ {p1}), \\
  \dot v &= \varepsilon ^{1-s} (b _ {j} u^ {j}+\sum_{k=1}^{j-1}b_ku^k),
  \qquad s=1/(j+1), \quad j\geq 2,
\end{split}
\end{equation}
where $b_k$ are the deformation parameters.  Systems in the form
(\ref{eq12}) were examined in \cite{Mor2002,Mor2004}.

When the perturbation is Hamiltonian and harmonic (contains only
one Fourier harmonics),  then function $\widetilde{A}(v,I_{pq})$
is also harmonic \cite{Mor2005}. System (\ref{eq94}) can be
reduced to a form with $\widetilde{A}(v,I_{pq})=a_{pq}\sin{pv}$.
We note that in this case functions $\widetilde P_0(v,I_{pq})$,
$Q_0(v,I_{pq})$ are also harmonic. Let us denote $\widetilde
P_0(v,I_{pq})=c_{pq}\sin{pv}$. Then, due to (\ref{8.28}), we have
$Q_0(v,I_{pq})=d_{pq}\cos{pv}$, and $c_{pq}=pd_{pq}$. For a
harmonic perturbation it was shown that $q=1$ \cite{Mor2005}.

Thus, for harmonic Hamiltonian perturbations, system (\ref{8.27})
is reduced to the form
\begin{equation}\label{8.30}
\begin{split}
\dot u &=\varepsilon ^{1-s}a_{p1}\sin{pv}+
\varepsilon c_{p1}u\sin{pv}  \\
 \dot v &=\varepsilon ^{1-s}b_{j}u^{j}+
 \varepsilon (b_{j+1}u^{j+1}+d_{p1}\cos{pv}),\, s=1/(1+j),
  \end{split}
\end{equation}
In (\ref{8.30}) we proceed with the slow time
$\tau=\varepsilon^{1-s}t$. The Hamiltonian $\overline{H}(u,v)$ of
system (\ref{8.30}) takes the form
\begin{equation}\label{eq17}
  \overline{H}(u,v)=\frac{b_ju^{j+1}}{j+1}+\frac{a_{p1}}{p}\cos{(pv)}+
  \varepsilon ^s\left(
  \frac{c_{p1}}{p}u\cos{(pv)}+\frac{b_{j+1}}{j+2}u^{j+2}\right), \quad j\geq 2.
\end{equation}
Deformations of the vector field are described by a system with
the Hamiltonian
\begin{equation}\label{eq18}
  \widetilde{H}(u,v)=\overline{H}(u,v) +\sum_{k=1}^{j-1}\frac{b_ku^{k+1}}{k+1}.
\end{equation}

We investigate further the averaged system
\begin{equation}\label{eq19}
  \frac{du}{d\tau }=-\frac{\partial \widetilde{H}}{\partial v},
  \quad \frac{dv}{d\tau }=\frac{\partial \widetilde{H}}{\partial u}
\end{equation}
for $j=2$. In \cite{Mor2008}, an example was considered where a
system  in the form (\ref{eq19}) with $j=2$ was derived. Some
phase portraits of system (\ref{eq19}) were presented in
\cite{Mor2008}, however the complete description was not found
there.

\subsection{Investigation of the averaged system. Bifurcations}

Let us change $pv\to v$ in system (\ref{eq19}) and denote
$a=a_{p1}$, $b=b_2$, $\mu_1=\varepsilon^{1/3}c_{p1}$, $\mu_2=b_1$. Then, assuming
 $b_3=0$, we come to the Hamiltonian system
\begin{equation}\label{m1}
\begin{split}
u^{\prime} &= a\sin{v}+\mu_1 u\sin{v}\\
 v^{\prime}&=p(bu^2 +\mu_2 u)+\mu_1 \cos{v}
 \end{split}
\end{equation}
with the Hamilton function $\widetilde
H(u,v)=p(\mu_2u^2/2+bu^3/3)+(a+\mu_1u)\cos{v}$.  Here the prime
denotes the derivative in $\tau$. To be definite, we set $a=2$,
$b=1$. Let us forget for a minute on the smallness of the
parameter $\varepsilon$.

System (\ref{m1}) has equilibria
$$O_1^{\pm}((-\mu_2\pm \sqrt{\mu_2^2-4\mu_1/p})/2,0),
$$
$$O_2^{\pm}((-\mu_2 \pm \sqrt{\mu_2^2+4\mu_1/p})/2,
\pi ).
$$
System (\ref{m1}) can have more equilibria
$$O_3\left(\frac{-2}{\mu_1},
\arccos{\frac{2}{\mu_1^2}\left(p\mu_2-\frac{2}{\mu_1}\right)}\right),
$$
$$O_4\left(\frac{-2}{\mu_1}, \pi -\arccos{\frac{2}{\mu_1^2}\left(p\mu_2-
\frac{2}{\mu_1}\right)}\right),
$$
if the following conditions are satisfied
\begin{equation}\label{m2}
-1\leq \frac{2}{\mu_1^2}\left(p\mu_2-\frac{2}{\mu_1}\right)\leq 1
.
\end{equation}

The characteristic equation has the form $\lambda^2+\Delta=0$
where $\Delta$ is the determinant of the Jacobi matrix, i.e.,
$\Delta=-\mu_1^2\sin^2{v_0}-(2u_0+p\mu_2)(2+\mu_1u_0)\cos{v_0}$,
here $u_0,\,v_0$ are the coordinates of the equilibrium. Hence,
system (\ref{m1}) can have simple equilibria of only two types:
saddles and centers. If $\Delta<0$ then the equilibrium is of the
saddle type and if $\Delta>0$, then it is a center.

Bifurcations of the equilibria $O_1^{\pm}$ occur when parameters cross the parabola
\begin{equation}\label{m3}
p\mu_2^2-4\mu_1 =0,
\end{equation}
and for the equilibria $O_2^{\pm}$ do on the parabola
\begin{equation}\label{m4}
p\mu_2^2+4\mu_1 =0.
\end{equation}
The equilibrium  become double (degenerate) with $\lambda_{1,2}=0$
at the bifurcation. Before the bifurcation two simple equilibria
exist: a saddle and a center corresponding to the same value of
$v$ and different values of $u$. We shall call such bifurcations
to be ``vertical''. Before the bifurcation, meander type curves
exist, in accordance to the terminology in \cite{Simo}.

Bifurcations of equilibria $O_3$, $O_4$ take place on the bifurcation curves
\begin{equation}\label{m5}
p\mu_2 =\frac{2}{\mu_1}\pm \frac{\mu_1^2}{2}.
\end{equation}
Here at the moment of bifurcation, one has a triple saddle
\cite{ALGM-1}. Unlike the previous case, the bifurcation here is
``horizontal'' (i.e. before the bifurcation, meander curves are absent,
see Fig.2, VIII). The triple saddle breaks into two simple saddles
and a center.

Besides the aforementioned local bifurcations, global bifurcations
of recon\-nec\-ting separatrices take place in system (\ref{m1}).
They  can may happen when the system has equilibria of the saddle
type lying on different levels $u=u_1$ and $u=u_2$. Let $v_1$ and
$v_2$ be corresponding coordinates $v$ for these saddles. The
equation the separatrix curves has the form
$\widetilde{H}(u,v)=h_1,$ for the saddle $(u_1,v_1)$ with
$h_1=\widetilde{H}(u_1,v_1)$ and the form $\widetilde{H}(u,v)=h_2$
for the saddle $(u_2,v_2)$ with $h_2=\widetilde H(u_2,v_2)$. The
moment of reconnection bifurcation is defined by the equality
\begin{equation}\label{m6}
h_1(\mu_1,\mu_2)=h_2(\mu_1,\mu_2) .
\end{equation}
The related bifurcation curve separates domains I  and II,
IV and VII, V and VI on Fig.1 ($\mu_2>0$). All bifurcation curves
(\ref{m3})-(\ref{m6}) are shown on Fig.1. The bifurcation curves
divide the plane of the parameters $(\mu_1,\mu_2)$ into 20
domains (for bounded $\mu_1$, $\mu_2$). Due to the symmetry of
these curves, it is sufficient to plot phase portraits from
12 domains in the upper half-plane $\mu_2>0$ (see Fig.1). Then the
following bifurcation scenarios can occur.

\begin{figure}[bth]\label{fig1}
\begin{center}
\includegraphics[scale=0.7]{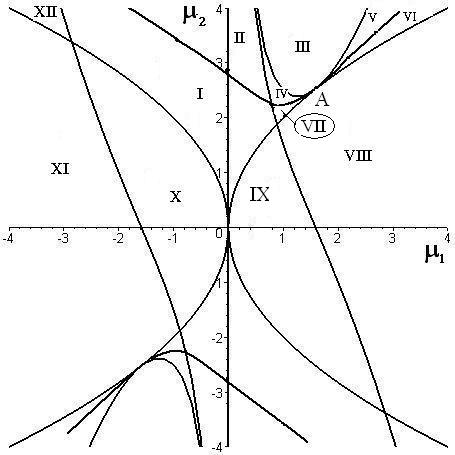}
\end{center}
\caption{The domains in the parameter plane
$(\mu_1,\mu_2)$ with different topologies of phase
portraits for system (\ref{m1}) at $a=2$, $b=1$ and $p=1$. }
\end{figure}

1) The first scenario is related to a transition from domain  I to
domain II (bifurcations of ``loops'').

2) The second scenario is related to formation (or annihilation)  of
``vortex pairs'' (a transition from domain IX to domains VIII).

3) The third scenario of codimension 2 is related to a transition
from domain III to domain VIII via the point A. At the bifurcation point,
there is a degenerate saddle with 6 separatrices
(\cite{ALGM-1}, pp. 385-404).

\begin{figure}[hp]\label{.
fig2}
\begin{center}
\begin{tabular}{ccc}
\footnotesize{(I)} & \footnotesize{(II)} & \footnotesize{(III)}\\
\includegraphics[scale=0.6]{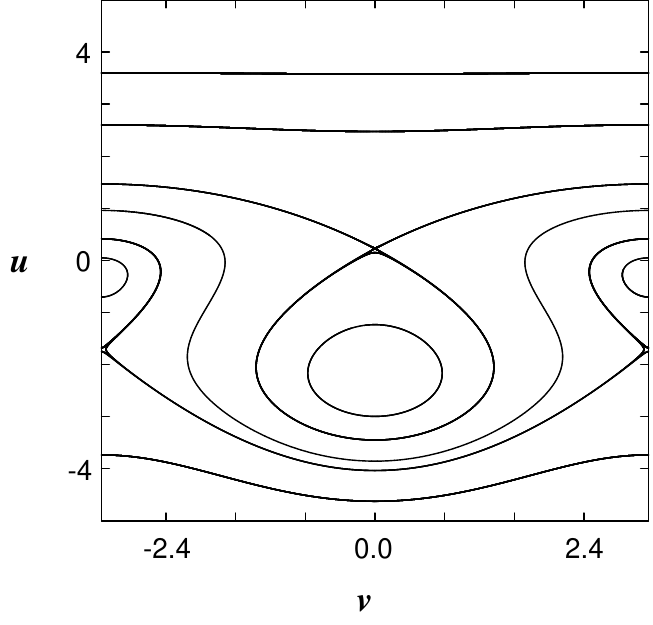}&
\includegraphics[scale=0.6]{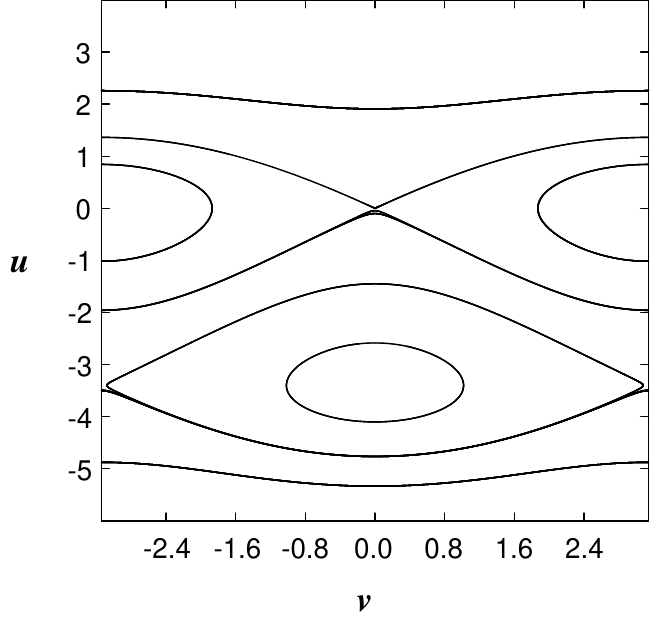}&
\includegraphics[scale=0.6]{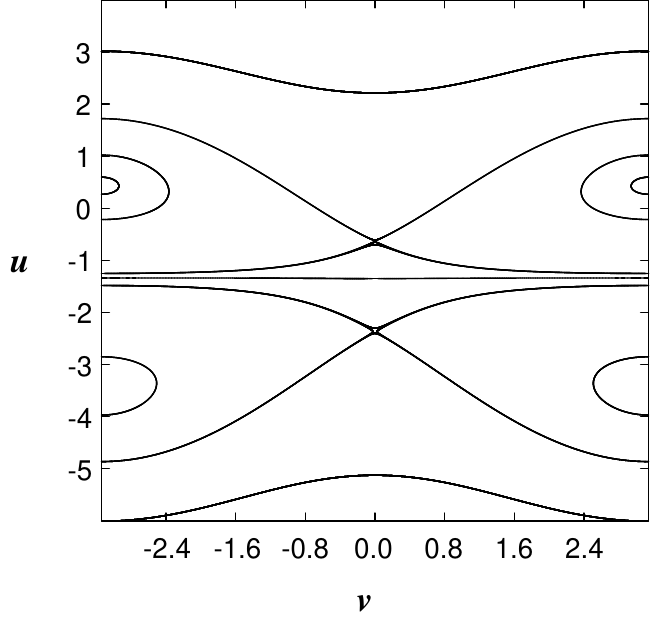}\\

\footnotesize{(IV)} & \footnotesize{(V)} & \footnotesize{(VI)}\\
\includegraphics[scale=0.6]{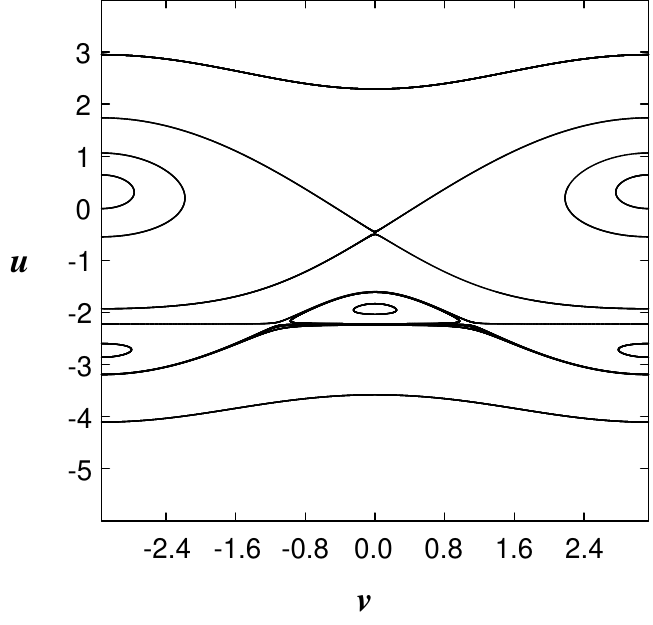}&
\includegraphics[scale=0.6]{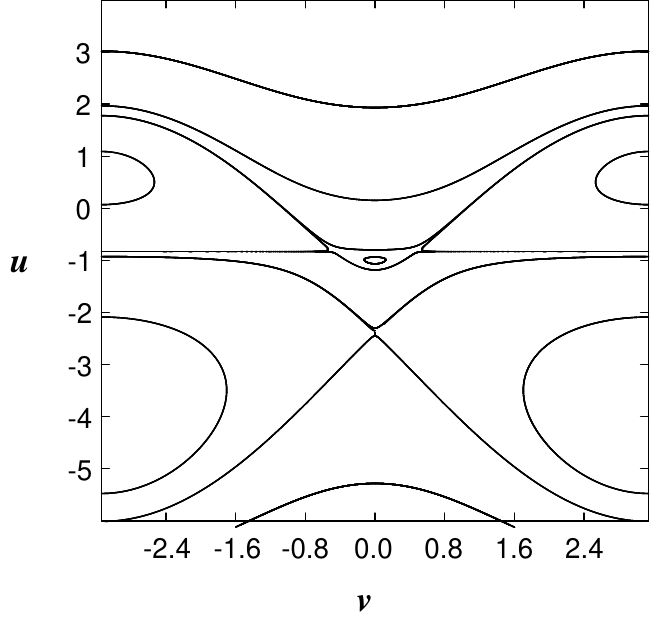}&
\includegraphics[scale=0.6]{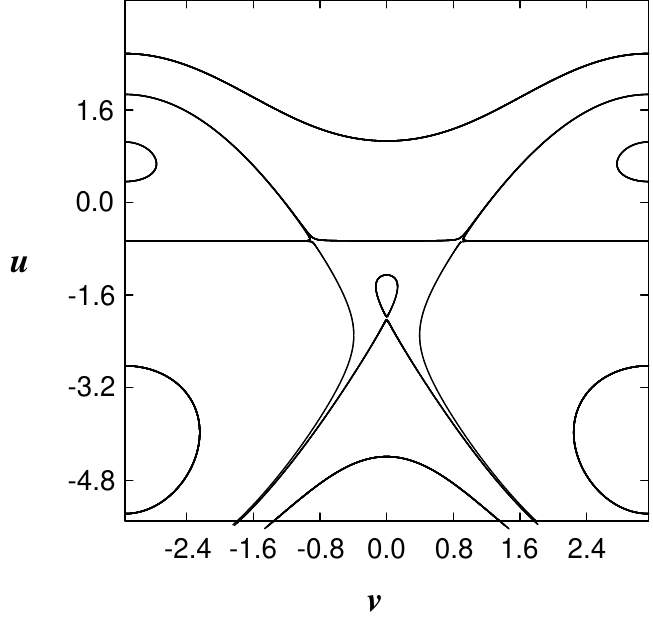}\\

\footnotesize{(VII)} & \footnotesize{(VIII)} & \footnotesize{(IX)}\\
\includegraphics[scale=0.6]{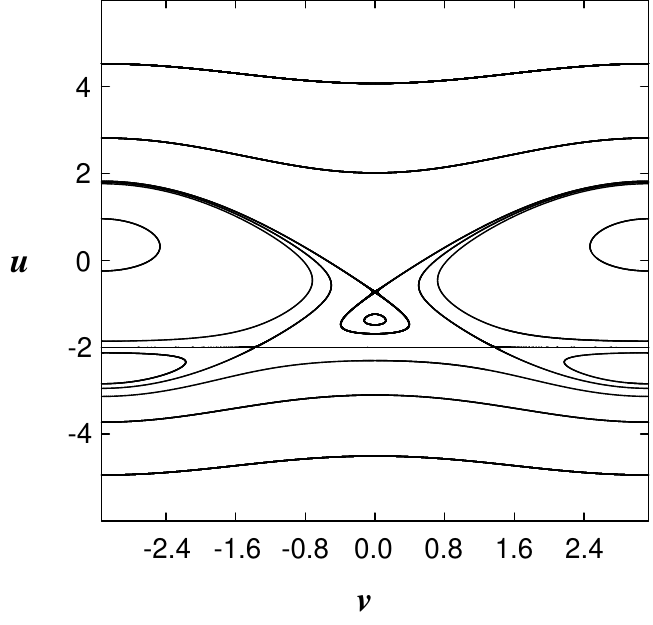}&
\includegraphics[scale=0.6]{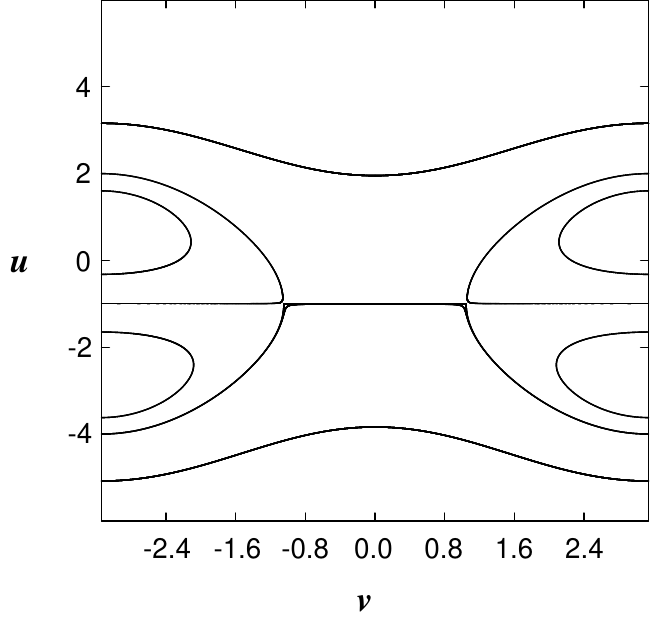}&
\includegraphics[scale=0.6]{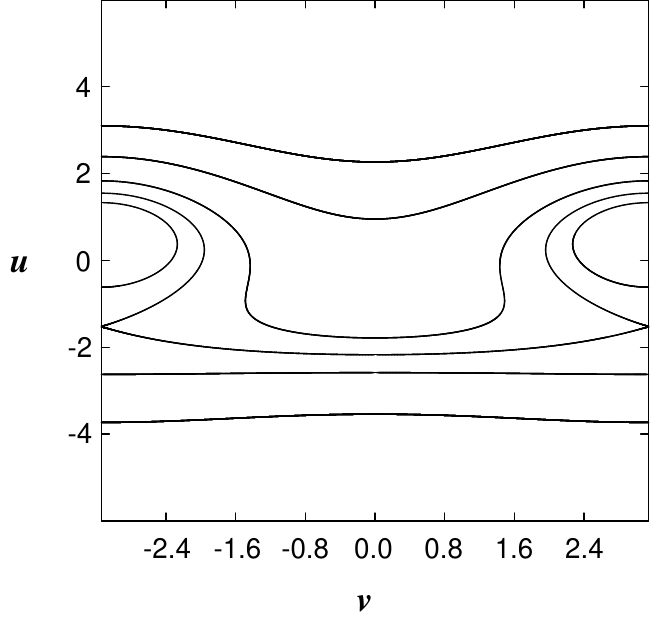}\\

\footnotesize{(X)} & \footnotesize{(XI)} & \footnotesize{(XII)}\\
\includegraphics[scale=0.6]{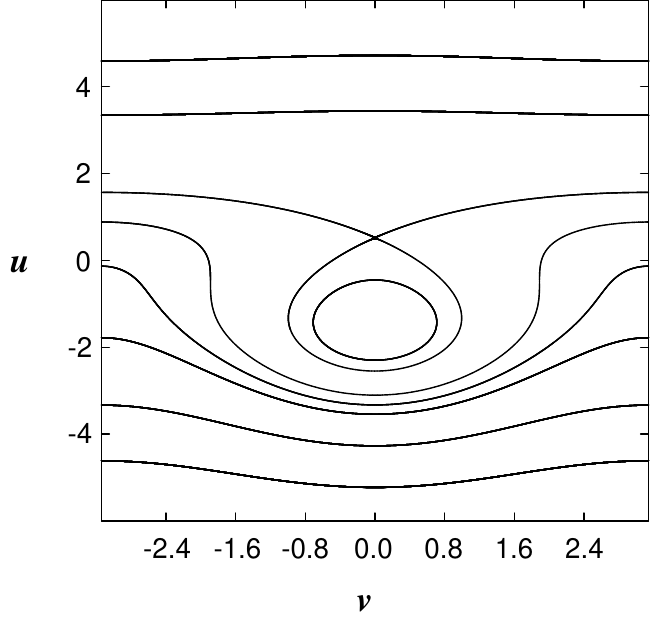}&
\includegraphics[scale=0.6]{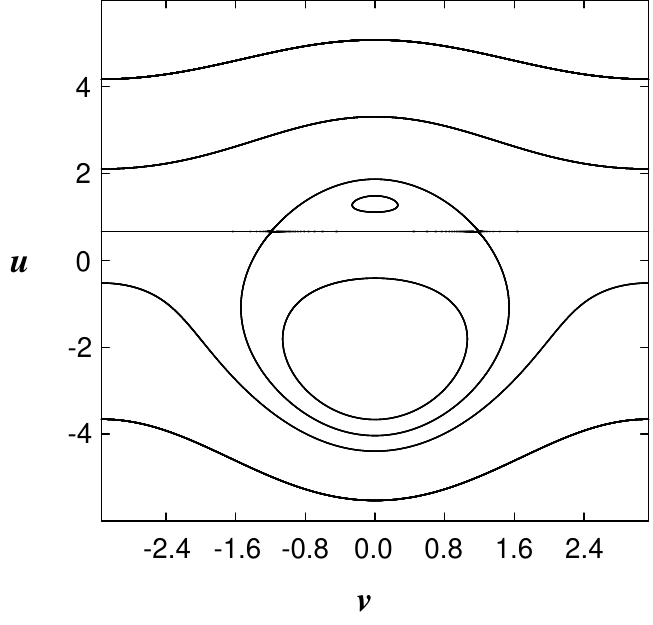}&
\includegraphics[scale=0.6]{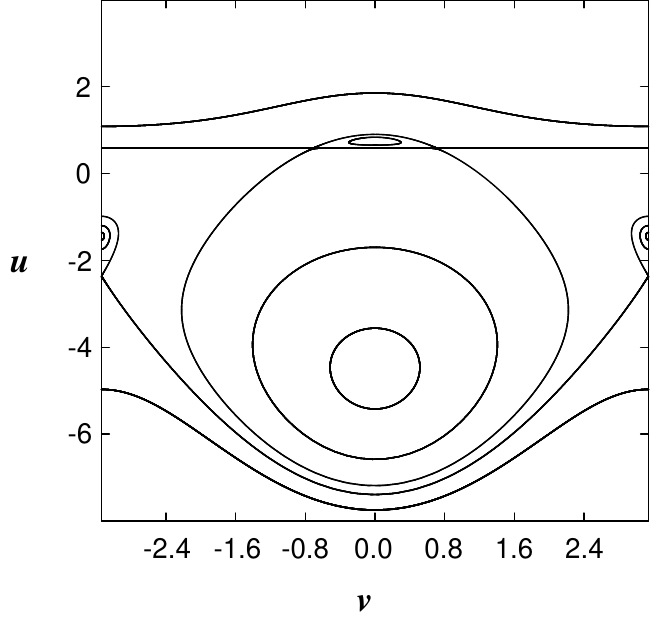}
\end{tabular}
\end{center}
\caption{Phase portraits of system (\ref{m1}).}
\end{figure}

Two types of reconnections are distinguished: 1) merging loops for
saddles  with $v=0$ and $u=\pi$ (for example, a transition from
domain I to domain II); 2) formation of a triangle from
separatrices of three saddles with coordinates $|v| <\pi$ (for
example, a transition  from domain IV to domain VII).

Certainly, other scenaria may exist. As an example,  we point
out the closed path on the bifurcation diagram related to a
transition from one domain to a neighboring one: $I\to II\to
IV\to III\to V\to VI\to VIII\to IX\to X\to XI\to XII\to I$ and
$I\to IX, I\to X; I\to VII \to VIII; IV \to VII$.  Bifurcations
that occur here can be understood using Fig.2.

 The bifurcation diagram on Fig.1 depends on the
parameters $a=a_{p1}$ and $b=d^2\omega(I_{p1})/2dI^2$. According
to \cite{Mor2005}, in the analytic case,
coefficient $a_{p1}$ exponentially decays with growing $p$. Hence,
for resonances with $p$ large, the transition from domain
I to domain VIII occurs for smaller values of
$\mu_1=\varepsilon^{1/3}c_{p1}$. The parameter $c_{p1}$ also
exponentially decays. Hence, trying to understand which scenario is realized, we
conclude that this depends on interrelations of
parameters $a$ and $\mu_1$ for the chosen parameters $b$ and
$\mu_2$. In any case, for $\mu_1$ small enough only the
first scenario is realized in system (\ref{m1}).

\section{Area preserving maps of the cylinder}

The averaged Hamiltonian system (\ref{m1}) is similar to systems
derived when analysing symplectic maps of the cylinder with a
non-monotone rotation \cite{HowH1984}-\cite{Petrisor}. In those
papers, two bifurcation scenarios were also found: ``loops'' and
``vortex pairs''. It is ascertained that the first scenario takes
place for the principle resonance $(p=1)$, while the second
scenario does for the sub-resonance of order 2 $(p=2)$. In
\cite{HowMor}, another map was considered being a combination of
the Chirikov and the Fermi maps. For this map, vortex pairs exist
also only for the sub-resonance of order 2.

In \cite{HowH1984,HowHum1995}, the standard map $T$ with a
non-monotone rotation depending on two parameters $a,\,\beta$ is
investigated
\begin{equation}\label{mm3}
\begin{split}
 u_{n+1}&=u_n + a\sin{v_n},\\
v_{n+1}&=v_n +u_{n+1}-\beta u_{n+1}^2,\quad n=0,1,2,...
\end{split}
\end{equation}
For this map,  the
approximating Hamilto\-ni\-an flow is generated by the Hamilton
function $H(u,v)=u^2/2-\beta u^3/3+a\cos{v}$, therefore only the first
scenario of ``loops'' is realized. For the map $T^2$, the
approximating Hamiltonian flow is generated by the Hamilton function
$H(u,v)=u^2/2-\beta u^3/3+(a^2/16)(1-2\beta u)\cos{2v}$, and the
second scenario of ``vortex pairs'' is realized. This bifurcation
is related to a transition from domain III to domain VIII through the
point $A$ (see Fig.1). Both Hamiltonians are special cases of
Hamiltonian (\ref{eq18}).

Following \cite{Mor2002}, consider the map
\begin{equation}\label{mm1}
\begin{split}
u_{n+1}&=u_n +\alpha (a+\mu_1u_{n+1})\sin{v_n}, \\
 v_{n+1}&=v_n
+\alpha (pbu_{n+1}^2+\mu_2pu_{n+1}+\mu_1 \cos{v_n}),\quad
n=0,1,2,...,
\end{split}
\end{equation}
that is derived by means of the conservative Euler method to solve
system (\ref{m1}).  Here $\alpha$ is the approximation step when the passage from
the differential system to the difference one, the error of the
approximation is of the order $\alpha^2$. We recast (\ref{mm1}) in
the form\footnote{Similar maps were considered in \cite{Simo}.}
\begin{equation}\label{mm2}
\begin{split}
 u_{n+1}&=(u_n +\alpha a\sin{v_n})/(1-\alpha \mu_1\sin{v_n}),\\
v_{n+1}&=v_n +\alpha (pbu_{n+1}^2+\mu_2 pu_{n+1}+\mu_1
\cos{v_n}),\quad n=0,1,2,... \end{split}
\end{equation}
For small $\alpha$, invariant curves of map (\ref{mm2}) are close
to orbits for system (\ref{m1}). Hence, the first and the second
bifurcation scenaria are realized. However, it it worth noting two
discrepancies: 1) if system (\ref{m1}) has a separatrix going from
a saddle to a saddle then the separatrices of map (\ref{mm2}) are
split generating a Poincar\'e homoclinic structure (Fig. 3a); 2)
the map has a boundary invariant curve embracing the cylinder (see
\cite{AWM} in this connection). Above this curve (on the upper
half-cylinder) or below it (on the lower half-cylinder), orbits
are wandering \cite{Mor2002} (see Fig. 3b).

\begin{figure}[bth]\label{fig4}
\begin{center}
\begin{tabular}{cc}
\includegraphics[scale=0.5]{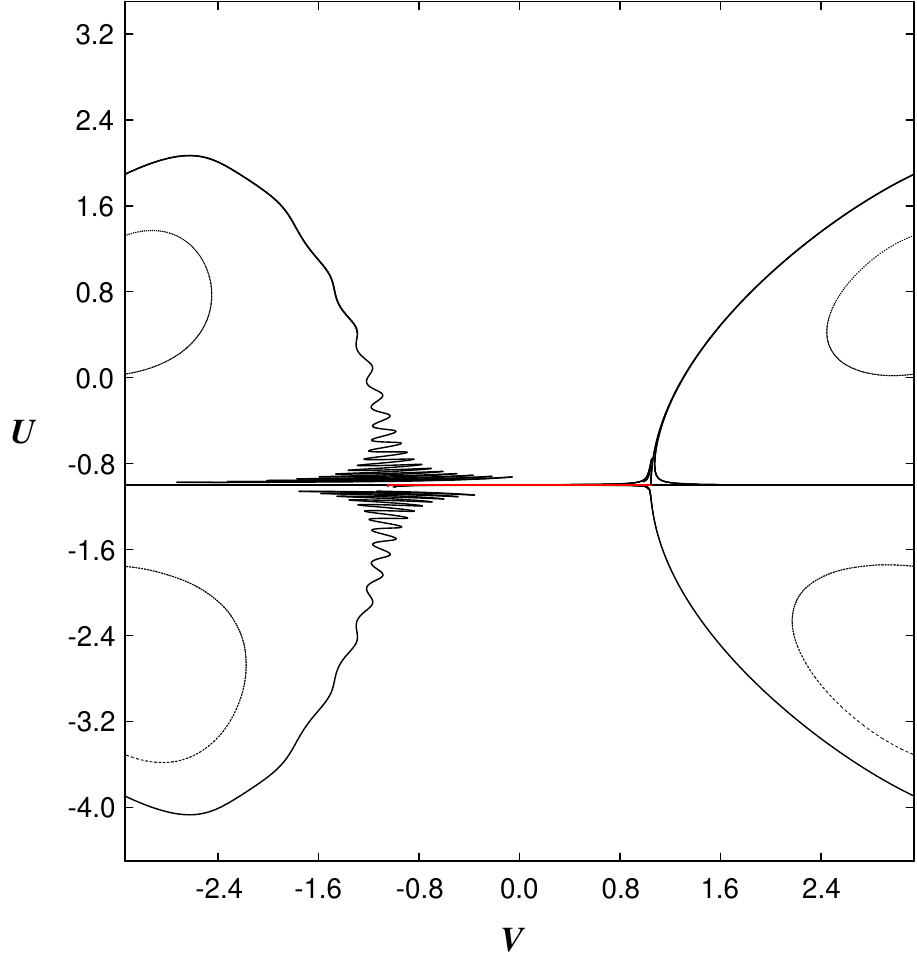}&
\includegraphics[scale=0.5]{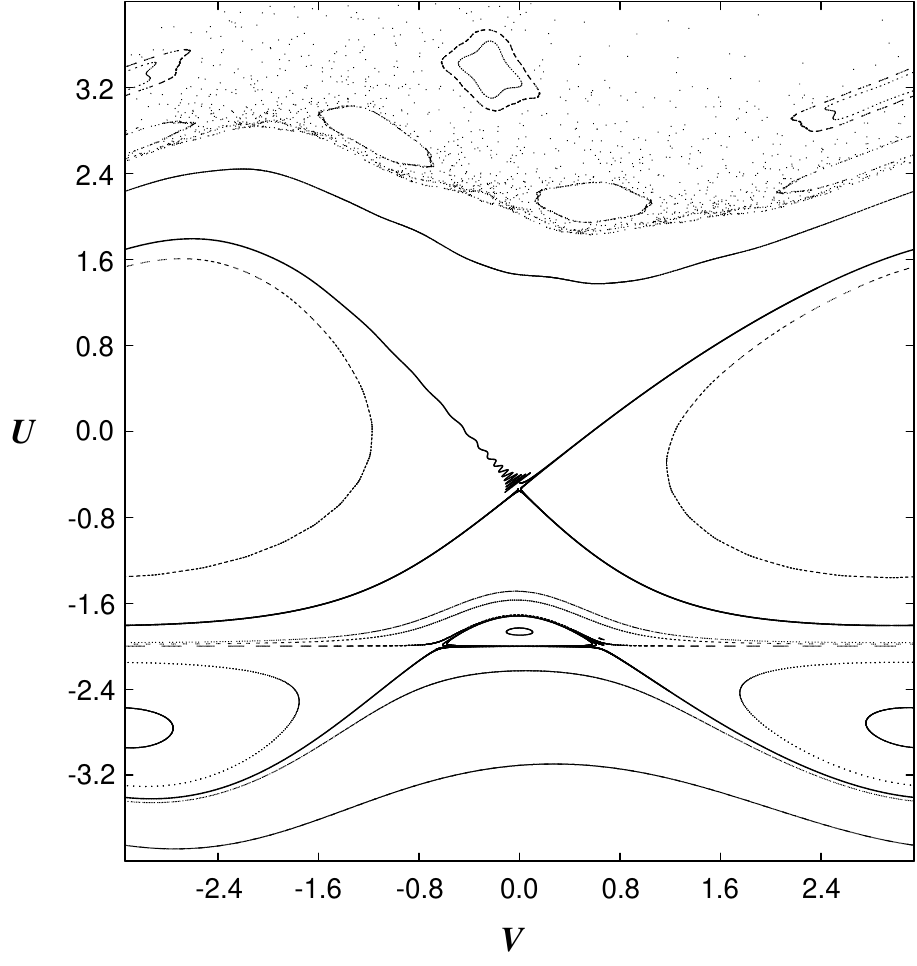}\\
\footnotesize{(a)} & \footnotesize{(b)}
\end{tabular}
\end{center}
\caption{The trajectories of map (\ref{mm2}) for $\alpha=0.17$,
$a=2$, $b=1$  corresponding to Fig.2 (VIII, IV). }
\end{figure}

\section{Conclusion}
The analysis of systems with 3/2 degrees of freedom being close to
two-dimensional integrable leads to the analysis of
two-dimensional systems (\ref{8.27}) on the cylinder. The latter
defines an approximating flow for the Poincar\'e map induced by
the initial system. In the case under consideration, the constructed map
(\ref{mm2}) defines the main features of the map.

For Hamiltonian systems studied above,  phase portraits
of system (\ref{m1}) provide all possible phase portraits
constructed for Hamiltonian flows approximating area preserving
maps of the cylinder with a non-monotone rotation
\cite{HowH1984}-\cite{Petrisor}. System (\ref{m1}) is a special
case of system (\ref{8.27}).

The work was partially supported by the Russian Science
Foundation, grant 14-41-00044.

\end{document}